\documentclass{ifacconf}
\usepackage{natbib,amsfonts,amsmath,graphicx}

\newcommand{\dotex}{{\frac{d}{dt}}}
\newcommand{\ddotex}{{\frac{d^2}{dt^2}}}

\begin{document}
\bibliographystyle{alpha}
\begin{frontmatter}

\title{Flatness-based  pre-compensation  \\
of  laser diodes .\thanksref{footnoteinfo}}

\thanks[footnoteinfo]{This work has been supported by the projet INRIA  STIC-Tunisie (2006/2007).   }

\author[First]{A. Abichou}
\author[Second]{S. El Asmi}
\author[Third]{P. Rouchon }

\address[First]{Laboratoire d'ing\'enierie math\'ematique (LIM), Ecole Polytechnique, BP 743, La Marsa, Tunisie (e-mail:
azgal.abichou@ept.rnu.tn)}
\address[Second]{SupCom Tunis, Cit\'e Technologique des Communications, 2083 Ariana (e-mail: elasmi@supcom.rnu.tn)}
\address[Third]{Ecole des Mines Paris, Centre Automatique et Syst\`{e}mes, 60, Bd-Saint-Michel, 75272 cedex 06, France (e-mail: pierre.rouchon@ensmp.fr)}


\begin{abstract}

A physical nonlinear dynamical model of a  laser diode is considered. We propose a
feed-forward control scheme   based on  differential flatness for the design of
 input-current modulations to compensate  diode distortions. The goal is to transform  without distortion  a radio-frequency current modulation  into a  light  modulation leaving the laser-diode and entering an optic fiber.  We  prove that  standard physical dynamical models based on dynamical electron and photons balance  are flat systems when the  current is considered as control input, the flat output being the photon number (proportional to the light power). We prove that input-current is an affine map of the flat output, its logarithm and their time-derivatives up to order two. When the flat output is an almost harmonic signal with slowly varying amplitude and phase, these derivatives admit precise analytic approximations. It is then possible to design   simple analogue  electronic circuits to code   approximations of  the  nonlinear computations required by our flatness-based approach.
Simulations with the parameters of a commercial diode illustrate the practical interest of this pre-compensation  scheme and its robustness versus modelling and analogue implementation  errors.

\end{abstract}

\end{frontmatter}

\section{Introduction}

Radio-over-Fiber is a very attractive technique for wireless access
network infrastructure, because it can transmit microwaves and
millimeter-waves through optical fibers for a long distance.
Therefore, the distribution of radio signals over optical fiber,
which  is of great interest for many applications such as broad-band
wireless access networks, sensor networks, radar and satellite
communications, has been intensively studied.  And there have been
rapid advances in the techniques to generate and transport radio
signals over optical fiber in recent years (see, e.g.,~\cite{Ai-Raweshidy-et-al:book02}). When the wireless channel is in series with the
optical link, nonlinear distortion due to the electrical/optical
conversion is the biggest concern (see, e.g.,~\cite{maury-et-al:ieee97}). This makes nonlinearity compensation an attractive
solution to improve link performance.

The focus of this work is to investigate the standard physical model used for describing commercial laser diodes. This model is  based on rate equation for the electron and photon populations. This allows us
 to propose a pre-compensation based on  differential flatness for the
design of the input current in order to cancel  diode distortions. Simulations corresponding to a commercial laser diode illustrate the pre-compensation scheme. When the distortion due to the
optic fiber is negligible, our  flatness-based pre-compensation scheme can be interesting to increase the bandwidth  and transmission rate.

This paper is organized as follows. In section~\ref{flat:sec} we recall the standard diode model and show that it is structurally flat, with the photon number being the flat output $y$. We show via elementary computations that the input current is a linear combination of $y$, $\dotex y$, $\ddotex y$, $\dotex(\log y)$ and $\ddotex (\log y)$ with coefficients depending via explicit formula on the physical parameters (equation~\eqref{flat:eq}). This relation is the starting point of the flatness based pre-compensation scheme addressed in section~\ref{precomp:sec}.  In Section~\ref{approx:sec}, we  propose  natural approximations when $y$ is quasi-harmonic  with slowly varying amplitude and phase.  These approximations  simplify considerably the computation burden in such a way that a nonlinear analogue circuit can be easily design to  realize  our flatness-based pre-compensation scheme.
In conclusion, we propose some hint to deal with phase influence  that becomes important when  fiber distortions cannot be negligible.

\section{The flat physical nonlinear model} \label{flat:sec}

A standard nonlinear model  relating the input-current $I$ to the
power of the emitted light  $P$ is recalled in \cite{Lee-et-al:ieeeQE03} where this model is used to analyze the nonlinear "transfer" between input-current modulation $I$ and the resulting output $P$) (see also \cite{agrawal-book1} for a simple exposure on such  basic models):
\begin{equation}\label{dynPN:eq}
 \left\{
 \begin{aligned}
   \dotex P & = R(P,X) - \frac{P}{\tau_p}
   \\
   \dotex X & = \frac{I}{I_{th} \tau_n} - \frac{X}{\tau_n} - a R(P,X)
 \end{aligned}
 \right.
\end{equation}
where
\begin{itemize}

 \item $\tau_p$ and $\tau_n$ are photon and carrier lifetimes;
 \item $X$ is the  normalized carrier  density;
 \item $I_{th}$ is the threshold current and $a$ is the constant related to other  physical parameters, $a=F I_{th} \tau_n/\tau_p$.

 \item $R(P,X)$ is the net rate of stimulated emission. Following~\cite{Lee-et-al:ieeeQE03} we use
 $$
 R(P,X)= \left(\frac{B \tau_n I_{th} (X-1) + \frac{1}{\tau_p}}{1+F B\tau_p\tau_c P}\right) P
 $$
 where $B$, $\tau_n$,$\tau_c$ are physical parameters .
\end{itemize}
Since the input-current appears only in the second differential equation of~\eqref{dynPN:eq} (the electron balance equation), the first state, $P$, is a flat output.
More precisely, this means that, if we consider~\eqref{dynPN:eq} as a control system with input $I$ and
output $y=P$, we  have  a nonlinear input/output  system  such that its inverse admits no dynamics. This system is differentially flat (see~\cite{fliess-et-al-cras92,fliess-et-al-ijc95,martin-et-al-caltech03,sira-agarwal:book04}) with $P$ as flat output:  the input $I$  is a
nonlinear function of $P$, $\dotex P$, $\ddotex P$. Notice that this
property is independent of the precise form of  $R(P,X)$ when it  depends effectively on $X$. This is always the case for physical reasons.

Let us detail now the computations with the precise form of $R$ given here above. The first equation of~\eqref{dynPN:eq}   yields
$$
 X=1+ \frac{P}{P_l} + \tau_p \dotex \left(\frac{ P}{P_l}\right)+  \tau_l \dotex \log\left(\frac{P}{P_l}\right)
$$
where $P_l= \frac{I_{th} \tau_n}{F \tau_c}$ and $\tau_l =\frac{1}{ B\tau_n I_{th}}$ are constant parameters.
Thus we can derive this expression versus $t$ to get $\dotex X$ as a
function of $P$, $\dotex P$ and $\ddotex P$:
$$
 \dotex X=\dotex\left(\frac{P}{P_l}\right) + \tau_p \ddotex \left(\frac{ P}{P_l}\right)+  \tau_l \ddotex \log\left(\frac{P}{P_l}\right)
$$
A linear combination  of the two
equations of~\eqref{dynPN:eq} gives $I$ explicitly:
$$
\frac{I}{I_{th}} = X + \tau_n \dot X  + \frac{F \tau_p}{\tau_n I_{th}} (\dot P + P/\tau_p)
.
$$
Thus we have the explicit formula relating $I$ to the  derivatives of $P$:
\begin{multline*}
      \frac{I}{I_{th}} =
       1 +
       \left(1+\frac{\tau_n}{\tau_c}\right) \left(\frac{P}{P_l}\right) +
      \left(\tau_n+\tau_p+\frac{\tau_n\tau_p}{\tau_c}\right) \dotex \left(\frac{P}{P_l}\right)
      \\ + \tau_n\tau_p \ddotex \left(\frac{P}{P_l}\right)+
      \tau_l \left(\dotex \log \left(\frac{P}{P_l}\right)  + \tau_n \ddotex \log \left(\frac{P}{P_l}\right)\right)
\end{multline*}
With the normalized input/ouput variables
\begin{equation}\label{uy:eq}
 u=\frac{I}{I_{th}}, \quad y = \frac{P}{P_l}= \frac{P}{\frac{I_{th} \tau_n}{F \tau_c}}
\end{equation}
 we have the following simple but nonlinear input/output relationships
that is equivalent to the physical model~\eqref{dynPN:eq}:
\begin{multline}\label{flat:eq}
  u =
       1 +
       \left(1+\frac{\tau_n}{\tau_c}\right) y +
      \left(\tau_n+\tau_p+\frac{\tau_n\tau_p}{\tau_c}\right) \dotex y
      \\ + \tau_n\tau_p \ddotex y +
      \tau_l \left(\dotex (\log y) + \tau_n \ddotex (\log y)\right)
\end{multline}
where the four parameters $\tau_n$, $\tau_p$, $\tau_c$ and $\tau_l=\frac{1}{ B\tau_n I_{th}}$ are positive time-scales. Remember that  the normalized electron density $X$ is a combination of $y$ and $\dotex y$:
$$
X=1+ y + \tau_p \dotex y+  \tau_l \dotex (\log y)
.
$$

\section{Flatness-based pre-compensation} \label{precomp:sec}

We can  now use~\eqref{flat:eq} for deriving feed-forward strategy  in order to compensate the non-linear distortion due to the diode.  Assume that our goal is to transfer the amplitude and phase modulated signal
$\epsilon(t) \cos(\omega t + \phi(t))$ from the current to the light leaving the diode. Here $\omega$ is the carriage pulsation, typically  around a few Ghz for radio frequencies, and $\epsilon(t)$ and $\phi(t)$ are slowly varying amplitude and phase where the bits are encoded. This means that
$\left|\dotex \epsilon\right| \ll \omega \epsilon $ and $\left|\dotex \phi\right| \ll \omega$. Without lost a generality, we can assume that $\epsilon < 1$. Assume finally that this modulated signal has to be converted in photons $y >0 $ according to $y = \bar y (1+\epsilon(t) \cos(\omega t + \phi(t)))$ where $\bar y$ is some positive constant.   This means that our goal is to find the current $u$ leading to such $y$. If one knows the diode parameters we can use~\eqref{flat:eq} to compute $u$.

\begin{figure}
$$  \includegraphics[width=.55\textwidth,height=.4\textheight]{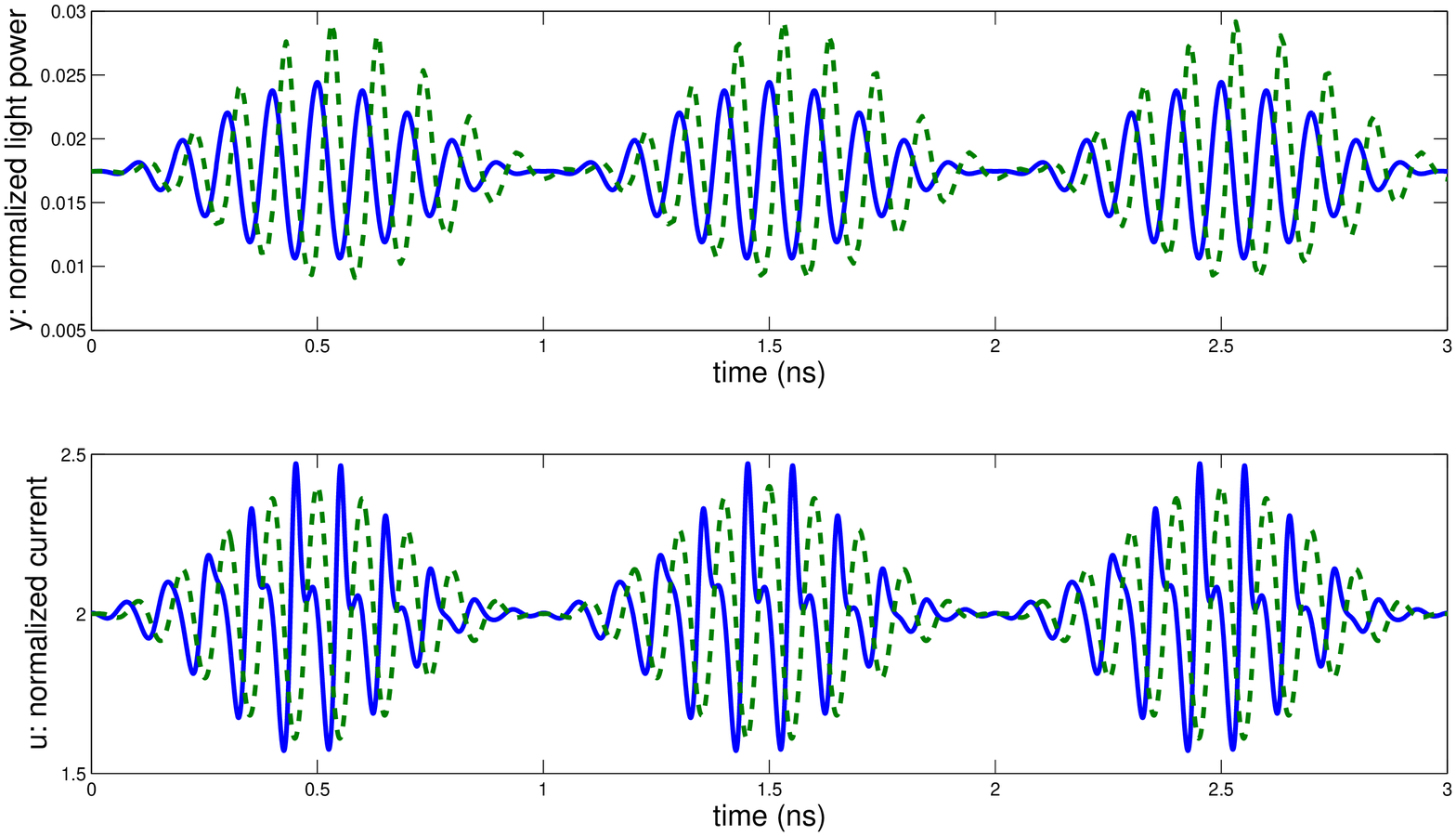}$$
  \caption{Input current modulations with $\epsilon(t)$ given by~\eqref{eps:eq}:
  the solid curve corresponds to the response of the diode when  $u$ is related to $y=\bar y (1+\epsilon(t) \cos(\omega t))$ via~\eqref{flat:eq}; dashed line corresponds to the response of the diode when
  $u= 1 + \left(1+\frac{\tau_n}{\tau_c}\right) \bar y (1+\epsilon(t) \cos(\omega t))$ (static model).   }\label{fig1}
\end{figure}

We have tested in simulation this simple pre-compensation scheme for the  commercial diode considered in~\cite{Lee-et-al:ieeeQE03}. For this diode,  we have
$$
\tau_n= 179~ps, \quad \tau_p =4.33 ~ps, \quad \tau_c=3.18~ps, \quad \tau_l=1.81~ps
.
$$
We choose a  modulation frequency $\omega$  of $10$~Ghz ($\omega=2\pi 10^{10}~s^{-1}$), $\phi=0$ and
\begin{equation}\label{eps:eq}
  \epsilon(t)= \frac{1}{5} \left(1-\cos\left(\frac{\omega t} {10}\right)\right)
  .
\end{equation}
The  average normalized light power is   $\bar y = 0.0175$ corresponds physically to $\bar P = 1.4$~mW since $P_l=78.5$~mW.  The simulations of figure~\ref{fig1}  illustrate the interest of taking into account the derivative terms in~\eqref{flat:eq} to modulate  the input  current $u$.

\section{Analogue implementation issues} \label{approx:sec}

For such high carrier  frequency $\omega \sim 10$~GHz, it is difficult and almost impossible  to use real-time numerical computations. Let us now propose two  approximations such that the above pre-compensation scheme can be applied in real-time via a specific electronic circuit.

\begin{figure}
$$  \includegraphics[width=.55\textwidth,height=.4\textheight]{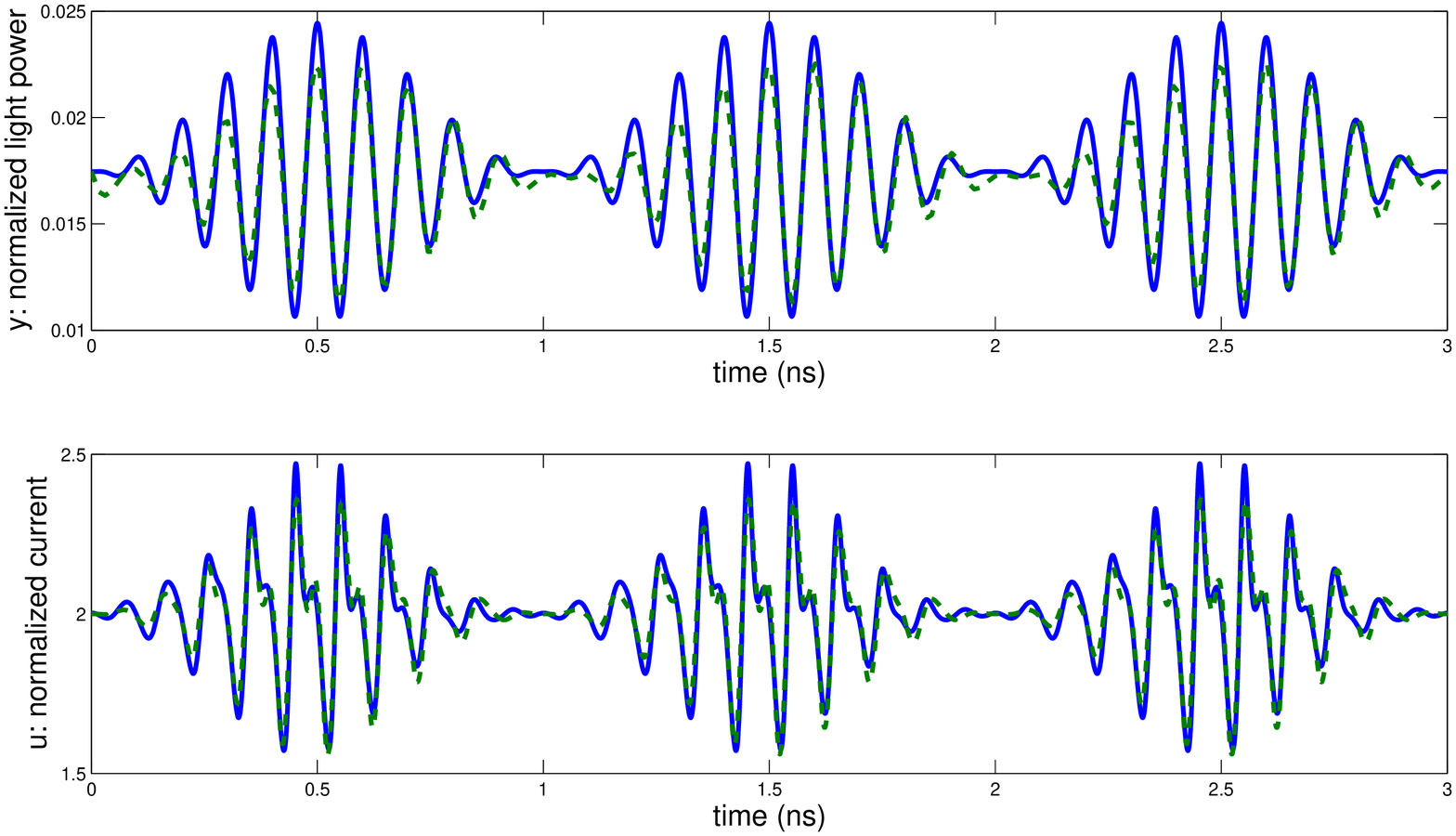}$$
  \caption{Input current modulations with $\epsilon(t)$ given by~\eqref{eps:eq}:
  the solid curve corresponds to the response of the diode when  $u$ is  perfectly related to $y=\bar y (1+\epsilon(t) \cos(\omega t))$ via~\eqref{flat:eq}; dashed line corresponds to the response of the diode when  approximations~\eqref{approx1:eq} and~\eqref{approx2:eq} are used in~\eqref{flat:eq}.}\label{fig2}
\end{figure}

\begin{figure}
$$  \includegraphics[width=.55\textwidth,height=.4\textheight]{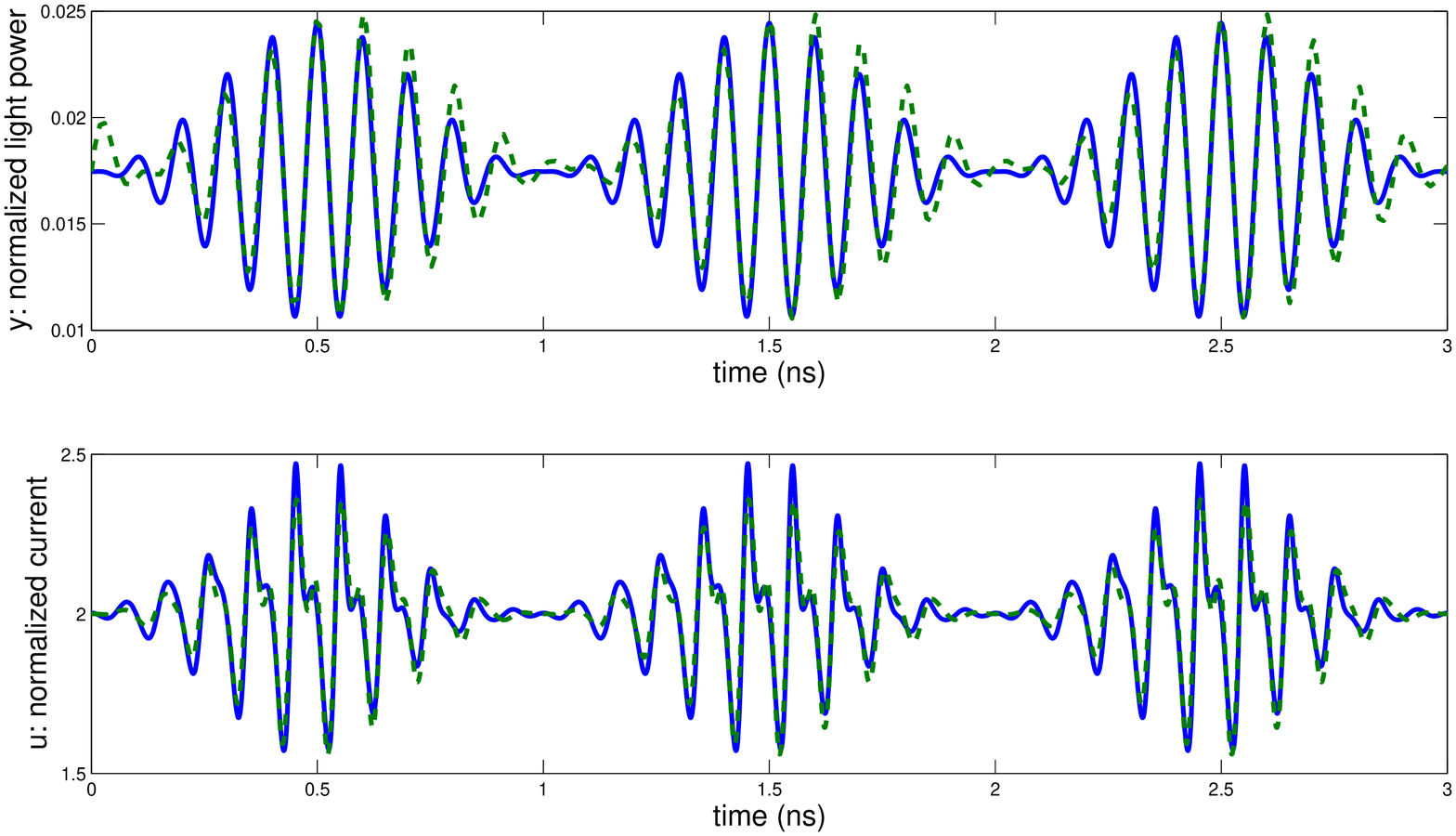}$$
  \caption{Input current modulations with $\epsilon(t)$ given by~\eqref{eps:eq}: as for figure~\ref{fig2},
  the solid curve corresponds to the response of the diode when  $u$ is  perfectly related to $y=\bar y (1+\epsilon(t) \cos(\omega t))$ via~\eqref{flat:eq}; dashed line corresponds to  approximations~\eqref{approx1:eq} and~\eqref{approx2:eq} plugged in~\eqref{flat:eq} but, contrarily to figure~\ref{fig2}, the carrier lifetime $\tau_n$ used for the computations of the diode response is multiplied by $2$ in order to test robustness versus large parametric errors.}\label{fig3}
\end{figure}

The first approximation is relative to the computation of derivative of $y$. Since we are looking for $u$ such that $y = \bar y (1+\epsilon(t)\cos(\omega t+\phi(t))$ with $\epsilon$ and $\phi$ slowly varying quantities, we have
\begin{align}\label{approx1:eq}
\dotex y \approx - \bar y \epsilon \omega \sin(\omega t +\phi), \qquad
\ddotex y \approx - \bar y \epsilon \omega^2 \cos(\omega t +\phi)
.
\end{align}
The second approximation is relative to the nonlinearity attached to the $\log$-terms in~\eqref{flat:eq}. Now we will use $\epsilon \ll 1$ and  propose an approximation up to second order terms:
\begin{align}
    \dotex (\log y) &\approx -\epsilon \omega \sin(\omega t + \phi)(1-\epsilon\cos(\omega t +\phi))
    \notag
    \\
    \ddotex (\log y) &\approx
    -\epsilon \omega^2 ( \cos(\omega t + \phi) + \epsilon (2\sin^2(\omega t + \phi)-1))
    \label{approx2:eq}
\end{align}
Plugging these approximations into~\eqref{flat:eq} shows that $u$ can be expressed approximatively  as a polynomial of degree $2$ in   $\epsilon \sin(\omega t +\phi)$ and $\epsilon  \cos(\omega t +\phi)$. Such  polynomial computations can be easily done by a nonlinear analogue circuit. Figure~\ref{fig2} shows that the diode response remains almost the same when such polynomial approximations are used. Other simulations displayed on figure~\ref{fig3} show that such pre-compensation scheme are also robust to  parameters uncertainties.

\section{Conclusion}

We have shown that it is possible to compute, via flatness-based motion planing  techniques,  adapted input current modulations in order to generate, in real-time and via analogue nonlinear circuits,  a given   modulation of the light  entering the optic fiber. If the distortion along the fiber is not negligible (long fiber and/or   refractive index  strongly dependent versus optic frequency), we have to consider also the phase dynamic of   the light emitted by the  diode (see, e.g.,~\cite{agrawal-book1} for more details about the spectrum widening  due to phase/amplitude variations).
Following~\cite{Lee-et-al:ieeeQE03} we have to complete~\eqref{dynPN:eq} by
$$
\dotex \varphi = \frac{\alpha B \tau_n I_th}{2} (X-1)
$$
where $\varphi$ is the phase of the emitted light and $\alpha$ a positive constant. Thus
the emitted light entering the optic fiber  is represented
by the  classical electric field (complex notations)
$$
    E(t)\propto \sqrt{P(t)} \exp(-\imath \varphi(t))\exp(-\imath\omega_0 t)
$$
where $\omega_0$ is the optical frequency (around $10^{15}$~Hz).
It is  interesting to notice that the diode dynamics (described by three states $(P,X,\varphi)$) and one input $I$)  is still flat with $z=\log P + F B \tau_p \tau_c P - \frac{2}{\alpha}\varphi$ being the  new flat output. This results from $\dotex z = - \frac{1 + F B \tau_p\tau_c P}{\tau_p}$. This structural property could certainly  be exploited in  other pre-compensation schemes taking into account  such  additional phase effect.


\begin{thebibliography}{8}
\providecommand{\natexlab}[1]{#1}
\providecommand{\url}[1]{\texttt{#1}}
\expandafter\ifx\csname urlstyle\endcsname\relax
  \providecommand{\doi}[1]{doi: #1}\else
  \providecommand{\doi}{doi: \begingroup \urlstyle{rm}\Url}\fi

\bibitem[Agrawal(1997)]{agrawal-book1}
G.P. Agrawal.
\newblock \emph{Fiber-Optic Communication Systems}.
\newblock John Wiley and Sons, 1997.

\bibitem[Ai-Raweshidy and Komaki(2002)]{Ai-Raweshidy-et-al:book02}
H.~Ai-Raweshidy and S.~Komaki.
\newblock \emph{Radio Over Fiber Technologies for Mobile Communications
  Networks}.
\newblock Artech House, Boston, 2002.

\bibitem[Fliess et~al.(1992)Fliess, L\'{e}vine, Martin, and
  Rouchon]{fliess-et-al-cras92}
M.~Fliess, J.~L\'{e}vine, Ph. Martin, and P.~Rouchon.
\newblock Sur les syst\`{e}mes non lin\'{e}aires diff\'{e}rentiellement plats.
\newblock \emph{C.R. Acad. Sci. Paris}, I--315:\penalty0 619--624, 1992.

\bibitem[Fliess et~al.(1995)Fliess, L\'{e}vine, Martin, and
  Rouchon]{fliess-et-al-ijc95}
M.~Fliess, J.~L\'{e}vine, Ph. Martin, and P.~Rouchon.
\newblock Flatness and defect of nonlinear systems: introductory theory and
  examples.
\newblock \emph{Int. J. Control}, 61\penalty0 (6):\penalty0 1327--1361, 1995.

\bibitem[Lee et~al.(2003)Lee, Choi, and Choi]{Lee-et-al:ieeeQE03}
K.-H. Lee, H.-Y. Choi, and W.-Y. Choi.
\newblock Analysis of chromatic dispersion-induced second-harmonic distortions
  including laser dynamics to the second order.
\newblock \emph{IEEE J. of Qauntum Electronics}, 39\penalty0 (5), 2003.

\bibitem[Martin et~al.(2003)Martin, Murray, and
  Rouchon]{martin-et-al-caltech03}
Ph. Martin, R.~Murray, and P.~Rouchon.
\newblock Flat systems, equivalence and trajectory generation, 2003.
\newblock Technical Report {\tt http://www.cds.caltech.edu/reports/}.

\bibitem[Maury et~al.({1997})Maury, Hilt, Berceli, Cabon, and
  Vilcot]{maury-et-al:ieee97}
G.~Maury, A.~Hilt, T.~Berceli, B.~Cabon, and A.~Vilcot.
\newblock Microwave-frequency conversion methods by optical interferometer and
  photodiode.
\newblock \emph{IEEE Trans. Microw. Theory Tech.}, {45}:\penalty0 {1481--1485},
  {1997}.

\bibitem[Sira-Ramirez and Agarwal(2004)]{sira-agarwal:book04}
H.~Sira-Ramirez and S.K. Agarwal.
\newblock \emph{Differentially flat systems}.
\newblock CRC, 2004.

\end{thebibliography}

\end{document}